\documentclass[12pt,twoside,reqno]{amsart}

\newtheorem{thm}{Theorem}[section]
\newtheorem{lemma}{Lemma}[section]
\newtheorem{remark}{Remark}[section]
\newtheorem{proposition}{Proposition}[section]
\newtheorem{corollary}{Corollary}[section]

\def \eps{\varepsilon}

\def \R{{\Bbb R}}

\numberwithin{equation}{section}

\begin{document}

\title[Subcritical and critical gZK]
{ Subcritical and critical generalized Zakharov-Kuznetsov equation
posed on bounded rectangles }
\author[
M. Castelli, \ G. Doronin] { M. Castelli$^\ast$, \ G. Doronin
%$^\S$
%$^\ast$
\bigskip
\\
{\tiny
Departamento de Matem\'atica,\\
Universidade Estadual de Maring\'a,\\
87020-900, Maring\'a - PR, Brazil. } }
\address
{
Departamento de Matem\'atica\\
Universidade Estadual de Maring\'a\\
87020-900, Maring\'a - PR, Brazil. }
\email{marcos\_castelli@hotmail.com \ \ ggdoronin@uem.br}
%\date{}

%\thanks{$^\ast$Supported by ....}
\thanks{$^\ast$Partially supported by CAPES}
%\thanks{$^\S$Corresponding author, partially supported by Funda\c{c}\~ao Arauc\'aria, Pr., Brazil}

\begin{abstract}
Initial-boundary value problem for the generalized
Zakharov-Kuznetsov equation posed on a bounded rectangle is
considered. Critical and subcritical powers in nonlinearity are
studied.
\end{abstract}

\maketitle

\section{Introduction}\label{introduction}

We are concerned with initial-boundary value problems (IBVPs) posed
on bounded rectangles located at the right half-plane
$\{(x,y)\in\mathbb{R}^2:\ x>0\}$  for the generalized
Zakharov-Kuznetsov \cite{pastor} equation
\begin{equation}
u_t+u_x+u^{1+\delta}u_x+u_{xxx}+u_{xyy}=0,\label{mzk}
\end{equation}
with $\delta \in [0,1].$ When $\delta=0,$ \eqref{mzk} turns the
classical Zakharov-Kuznetsov (ZK) equation \cite{zk}, while
$\delta=1$ corresponds to so-called modified Zakharov-Kuznetsov
(mZK) equation \cite{pastor1} which is a two-dimensional analog of
the well-known modified Korteweg-de Vries (mKdV) equation
\cite{bona2}
\begin{equation}\label{mkdv}
u_t+u_x+u^2u_x+u_{xxx}=0.
\end{equation}
Notes that both ZK and mZK possess real plasma physics applications
\cite{zk}.

As far as ZK is concerned, the results on both IVP and IBVPs can be
found in
%readers who are interested in both IVP and IBVP are referred to
\cite{faminski,faminski2,farah,pastor,pastor2,saut,temam, temam2}.
For IVP to mZK, see \cite{pastor1}; at the same time we do not know
solid results concerning IBVP to mZK. The main difference between
initial and initial-boundary value problems is that IVP provides
(almost immediately) good estimates in $(L^{\infty}_t;H^1_{xy})$ by
the conservation laws, while IBVP does not possesses this advantage.

Our work is a natural continuation of \cite{doronin} where
\eqref{mzk} with $\delta =0$ has been considered. There one can find
out a more detailed background, descriptions of main features and
the deployed reference list.

In the present note we put forward an analysis of \eqref{mzk} for
$\delta \in (0,1].$ When $\delta = 1,$ the power is critical (see
\cite{pastor,pastor1}) and a challenge concerning the well-posedness
of IBVPs appears. For one-dimensional dispersive models the critical
nonlinearity has been treated in \cite{larkin19}.

Once $\delta\in (0,1)$ the existence of a weak solution in
$\left((L^{\infty}_T; L^2)\cap (L^2_T;H^1_0)\right)$ with $u_0\in
L^2_{xy}$ is proved in our work via parabolic regularization. If
$\delta =1,$ we apply the fixed point arguments to prove the local
existence and uniqueness of solutions with more regular initial
data. We also show the exponential decay of $L^2$ norm of solutions
as $t\to\infty$ if $u\in (L^{\infty}_{\mathbb{R}^+} ; H^1_0),$ under
domain's size restrictions. These are the main results of the paper.

\section{Problem and notations}\label{problem}

%$\overset{c}{\hookrightarrow}$

Let $L,B,T$ be finite positive numbers. Define $\Omega$ and $Q_T$ to
be spatial and time-spatial domains
\begin{equation*}
\Omega=\{(x,y)\in\mathbb{R}^2: \ x\in(0,L),\ y\in(-B,B) \},\ \ \
Q_T=\Omega \times (0,T).
\end{equation*}

In $Q_T$ we consider the following IBVP:
\begin{align}
A&u\equiv u_t+u_x+u^{1+\delta}u_x+u_{xxx}+u_{xyy}=0,\ \ \text{in}\
Q_T; \label{2.1}
\\
&u(x,-B,t)=u(x,B,t)=0,\ \ x\in(0,L),\ t>0;
\label{2.2}
\\
&u(0,y,t)=u(L,y,t)=u_x(L,y,t)=0,\ \ y\in(-B,B),\ t>0;
\label{2.3}
\\
&u(x,y,0)=u_0(x,y),\ \ (x,y)\in\Omega, \label{2.4}
\end{align}
where $u_0:\Omega\to\mathbb{R}$ is a given function.

Hereafter subscripts $u_x,\ u_{xy},$ etc. denote the partial
derivatives, as well as $\partial_x$ or $\partial_{xy}^2$ when it is
convenient. Operators $\nabla$ and $\Delta$ are the gradient and
Laplacian acting over $\Omega.$ By $(\cdot,\cdot)$ and $\|\cdot\|$
we denote the inner product and the norm in $L^2(\Omega),$ and
$\|\cdot\|_{H^k}$ stands for the norm in $L^2$-based Sobolev spaces.
Abbreviations like $(L^s_t;L^l_{xy})$ are also used for anisotropic
spaces.

\section{Existence in sub-critical case}\label{existence}

In this section we state the existence result in sub-critical case,
i.e., for $\delta\in(0,1).$ We provide a short motivation for this
study at the final of the section.

\subsection{Sub-critical nonlinearity}
\begin{thm}\label{theorem1}
Let $\delta\in(0,1)$ and $u_0\in L^2(\Omega)$ be a given function.
Then for all finite positive $B,\ L,\ T$ there exists a weak
solution to \eqref{2.1}-\eqref{2.4} such that
\begin{equation*}
u\in L^{\infty}(0,T;L^2(\Omega))\cap L^2(0,T;H^1_0(\Omega)).
\end{equation*}
\end{thm}

To prove this theorem we consider for all real $\eps>0$
the following parabolic regularization of \eqref{2.1}-\eqref{2.4}:
\begin{align}
A^{\eps}u_{\eps}&\equiv Au_{\eps}+\eps(\partial_x^4u_{\eps}+\partial_y^4u_{\eps})=0\ \ \text{in}\ Q_T;\label{3.1}\\
&u_{\eps}(x,-B,t)=u_{\eps}(x,B,t)
% \notag\\&
=\partial_y^2u_{\eps}(x,-B,t)=\partial^2_yu_{\eps}(x,B,t)=0,\ x\in(0,L),\ t>0;\label{3.2}\\
&u_{\eps}(0,y,t)=u_{\eps}(L,y,t)
%\notag\\&
=\partial_x^2u_{\eps}(0,y,t)=\partial_xu_{\eps}(L,y,t)=0,\ y\in (-B,B),\ t>0;\label{3.3}\\
&u_{\eps}(x,y,0)=u_{0}(x,y),\ (x,y)\in \Omega.\label{3.4}
\end{align}

For all $\eps>0,$ \eqref{3.1}-\eqref{3.4} admits
%, at least for small $T>0$,
a unique regular solution in $Q_T$ \cite{lady2}. In what follows we
omit the subscript $\eps$ whenever it is unambiguous.

%\subsection*{Estimate I}

Multiplying $A^{\eps}u_{\eps}$ by $u_{\eps}$ and integrating over
$Q_T,$ we have

%\begin{equation}\label{estimate 1}
%\dfrac{d}{dt} \| u\| ^2 (t) + \int_{-B}^B u^2_x (0,y,t) \, dy +
%2\eps \big( \|u_{xx}\|^2 (t) + \|u_{yy}\|^2 (t) \big) = 0.
%\end{equation}
%Integrating in $t\in(0,T)$ gives

\begin{equation}\label{estimate 2.1}
 \| u\| ^2 (t) + \int_0^t \int_{-B}^B u^2_x (0,y,\tau) \, dyd\tau
 + 2\epsilon \int_0^t \big( \|u_{xx}\|^2 (\tau) + \|u_{yy}\|^2 (\tau) \big)d\tau = \|u_0\|^2, \,\,t\in(0,T).
\end{equation}

Multiplying $A^{\eps}u_{\eps}$ by $x u_{\eps}$, integrating over
$\Omega$ with the use of the Nirenberg, H\"older and Young
inequalities yields
\begin{eqnarray}\label{estimate 2.3}
\dfrac{d}{dt} \| \sqrt{x}u\|^2(t) +\frac{1}{2} \|\nabla u \|^2(t) +
2\|u_x\|^2(t) + 2\eps\Big( \|\sqrt{x}u_{xx}\|^2(t) +  \|\sqrt{x} u_{yy}\|^2 (t)\Big)  \nonumber \\
\leq \|u\|^2(t)+2\eps\int_{-B}^B u_x^2(0,y,t)\, dy  +
 \dfrac{C(\xi,\delta)C_{\Omega}^{\frac{2}{1-\delta}}}{3+\delta}\|u\|^{\frac{4}{1-\delta}}(t).
\end{eqnarray}
Integrating with respect to $t>0$ in \eqref{estimate 2.3} and taking
$\eps<1/2$ gives
\begin{eqnarray}\label{estimate 22.3}
&\| \sqrt{x}u\|^2(t)+\frac{1}{2}\int_0^{t} \|\nabla u \|^2(\tau)\,
d\tau + 2\int_0^{t} \|u_x\|^2(\tau)\, d\tau + 2\eps\int_0^{t}\Big(
\|\sqrt{x}u_{xx}\|^2(\tau) +  \|\sqrt{x} u_{yy}\|^2 (\tau)\Big)\,
d\tau
\nonumber \\
&\leq \int_0^{t}\|u_0\|^2\, d\tau +\int_0^{t}\int_{-B}^B
u_x^2(0,y,\tau)\,
dyd\tau+\dfrac{C(\xi,\delta)C_{\Omega}^{\frac{2}{1-\delta}}}{3+\delta}\cdot\int_0^{t}\|u_0\|^{\frac{4}{1-\delta}}d\tau
\nonumber\\
&\leq (T+1)\|u_0\|^2 +
\dfrac{C(\xi,\delta)C_{\Omega}^{\frac{2}{1-\delta}}}{3+\delta}\cdot
T  \|u_0\|^{\frac{4}{1-\delta}}.
%\nonumber \\
\end{eqnarray}
\begin{remark}
Note that \eqref{estimate 22.3} does not hold for critical
case, i.e., while $\delta \to 1.$
\end{remark}
Estimates \eqref{estimate 2.1} and \eqref{estimate 22.3} thus become
\begin{equation}\label{limitations}
\begin{array}{c}
\hspace{-2cm}u_{\eps} \;\;\;\; \text{is bounded in } \;\;\;\; L^{\infty}\big(0,T;L^2(\Omega)\big) ,\\
u_{\eps x}(0,y,t) \;\;\;\; \text{is bounded in } \;\;\;\; L^2\big(0,T;L^2(-B,B)\big) ,\\
\hspace{-2cm}\nabla u_{\eps}\;\;\;\; \text{is bounded in } \;\;\;\;
L^2\big(0,T;L^2(\Omega)\big),
\end{array}
\end{equation}
where limitations do not depend on $\eps$ but depend only on $T$,
$\delta$, $\Omega$ and $\|u_0\|$.

Thanks to (\ref{limitations}) we have boundness of
$u_{\eps}^{1+\delta}u_{\eps x}$ for all $\delta \in (0,1)$. In fact,
given $\delta \in (0,1)$ take $m=\frac{4}{3+\delta}$ and $ \kappa
(\delta) = \frac{1+\delta}{3+\delta}.$ Then H\"{o}lder's and
Nirenberg's inequality yield
\begin{eqnarray}\label{estimativa 3.3}
\|u^{1+\delta}u_x\|^m_{L^m(0,T;L^m(\Omega))} &= \int_0^T
\|u^{1+\delta}u_x\|^m_{L^m(\Omega)}(t)\, dt
%\nonumber \\&
\leq  C_{\Omega}^{4 \kappa (\delta)} \int_0^T \| \nabla u \|^{2 } (t)\|u \|^{2 \kappa (\delta)}(t) \, dt
\nonumber \\
& = C_{\Omega}^{4 \kappa (\delta)} \|u \|^{2 \kappa
(\delta)}_{L^{\infty}(0,T;L^2(\Omega))}  \| \nabla u \|^{2 }
_{L^2(0,T;L^2(\Omega))} .
\end{eqnarray}
Therefore, due to (\ref{estimativa 3.3}) and (\ref{limitations}) we
conclude that $u^{1+\delta}u_x$ is bounded in
$L^m(0,T;L^m(\Omega)).$ Since $L^{\frac{4}{1-\delta}}$ is the dual
space of $L^{ \frac{4}{3+\delta} } $ and $ H^{ 1 } \subset
L^{\frac{4}{1-\delta}}$ in dimension 2, we have as well

\begin{equation}\label{estimativa 3.5}
u^{1+\delta}u_x \;\;\;\text{is bounded in}\;\;\;
L^{\frac{4}{3+\delta}}(0,T;H^{-1}(\Omega)).
\end{equation}
Thanks to (\ref{limitations}) and (\ref{estimativa 3.5}) jointly
with the equation, we get
\begin{equation}\label{estimativa 3.6}
\dfrac{\partial u_{\epsilon}}{\partial t } \;\;\;\text{is bounded
(independently of $\eps$) in}
\;\;\;L^{\frac{4}{3+\delta}}(0,T;H^{-3}(\Omega))
\end{equation}
which assures the family $u_{\eps}$ to be relatively compact in
$L^2(0,T;L^2(\Omega))$. This is sufficiently to obtain the existence
of $\lim u_{\eps}$ as $\eps \rightarrow 0$, using the compactness
argument in the nonlinear term.

The initial condition $u(x,y,0) = u_0(x,y)$ is fulfilled; indeed,
due to (\ref{estimativa 3.6}) $ u_{\eps}$ converges to $u$ in
$C\big([0,T]; H^{-3}_{w}(\Omega)\big), $ where $H^{-3}_{w}$ is
$H^{-3}$ equipped with the weak topology.

By the same way, the Dirichlet condition $u=0$ onto $\partial
\Omega$ is satisfied since $u_{\eps}$ converges to $u$ weakly in
$L^{2}(0,T;H^{1}_0(\Omega)).$ It remains to show that
%\begin{equation}\label{cond cont}
$u_x(L,y,t)=0,$
%\end{equation}
which is
%can be done in the similar way as in
done by the following two lemmas (cf. \cite{temam, temam2}).

\begin{lemma}\label{lema cont} If $u\in  L^{\infty }(0,T;L^2(\Omega))\cap L^2 (0,T;H^1_0 (\Omega) )$
solves (\ref{2.1}), then
\begin{equation}\label{eq 1 lm cont}
u_x , u_{xx} \in C(0,L; V)\ \text{ with }\ V= H^{-2}\big((0,T)\times
(-B,B)\big),
\end{equation}
and, in particular,
\begin{equation}\label{eq 2 lm cont}
u_x \big|_{x=0,1}, \;\; u_{xx}\big|_{x=0,1}
\end{equation}
are well defined in $V$. Moreover, these traces depend continuously
of $u$ in an appropriate sense.
%Alem disso, esses tracos dependem
%continuamente de $u$ em um sentido preciso na prova.
\end{lemma}

To prove this lemma, write (\ref{2.1}) in the form
\begin{equation}\label{lm eq 1}
u_{xxx} = -u_x - u_{xyy} - u^{1+\delta}u_x - u_t ,
\end{equation}
and observe that
\begin{eqnarray*}
u_t \in L^{2}(0,L;H^{-1}\big(0,T;L^2 (-B,B)\big), \\
u_{xyy} \in L^{2}(0,L;L^2\big(0,T;H^{-2} (-B,B)\big).
\end{eqnarray*}
Accordingly with (\ref{estimativa 3.5}) and definition of $V$ in
\eqref{eq 1 lm cont}, it holds
\begin{equation}\label{lm eq 2}
u^{1+\delta}u_x \in
L^{\frac{4}{3+\delta}}\big(0,L;L^{\frac{4}{3+\delta}}((0,T)\times(-B,B))\big)
\hookrightarrow L^{\frac{4}{3+\delta}}\big(0,L;V\big).
%H^{-2}((0,T)\times(-B,B))\big).
\end{equation}
Thus we have
%$u_{xxx}$ pertence (pelo menos) ao maior desses
%espacos, isto eh,
\begin{equation}\label{lm eq 3}
u_{xxx} \in L^{\frac{4}{3+\delta}}\big(0,L;V\big)
%H^{-2}((0,T)\times(-B,B))\big),
\end{equation}
and (\ref{eq 1 lm cont}) and (\ref{eq 2 lm cont}) follow. Moreover,
if a sequence of functions $u_m$ satisfies (\ref{eq mZK}) and $u_m
\rightarrow u  $ in $ L^{\infty }(0,T;L^2(\Omega))\cap L^2
(0,T;H^1_0 (\Omega) )  $  strongly, then $ u_{mx} \big|_{x=0,1},
\;\; u_{mxx}\big|_{x=0,1}  $ converge to $u_x \big|_{x=0,1}, \;\;
u_{xx}\big|_{x=0,1}$ in $V.$ If a convergence of $u_m$ being weak
(star-weak for $L^{\infty }$,) then a convergence take place in
%sera mantida em
$C(0,L; V_w)$ and $Y_w$. This is based on compactness arguments
justified by (\ref{estimativa 3.6}), used to prove that
$u_m^{1+\delta}u_{mx} \rightarrow u^{1+\delta}u_x$.

\begin{lemma}\label{lema ux(L)}Let $U$ be a reflexive Banach space and $p\geq 1$.
Suppose that two function sequences $u_{\eps} ,\ g_{\eps}\in
L^p(0,L;U)$ satisfy
\begin{equation}\label{eq traco ux(L)}
 \left.
\begin{array}{c}
%-------------------------------
u_{\eps xxx} + \eps u_{\eps xxxx}= g_{\eps}, \\
%-------------------------------
u_{\eps}(0)=u_{\eps}(L)=u_{\eps x}(L)=u_{\eps xx}(0)=0 ,
%-------------------------------
\end{array}
\right.
\end{equation}
with $g_{\eps}$ being bounded in $L^p(0,L;U)$ as $\eps \to 0.$ Then
$ u_{\eps xx}  $ (consequently $ u_{\eps x}, $ and $ u_{\eps} $) is
bounded in $L^{\infty}(0,L;U)$ as $\eps \rightarrow 0.$ Moreover,
for a subsequence $u_{\eps} \to u $ converging (strongly or weakly)
in $L^q(0,L;U),$ $1\leq q < \infty,$ it holds that $u_{\eps x}(L)$
converges to $u_x(L)$ in $U$ (at least weakly), and therefore
$u_x(L)=0.$
\end{lemma}

See \cite{temam2} for the proof.

To prove Theorem \ref{theorem1}, apply the above lemmas with
$$
g_{\eps} := - u_{\eps t} - \eps u_{\eps x} - u_{\eps xyy} -
u^{1+\delta}_{\eps}u_{\eps \eps} - \eps u_{\eps yyyy} ,
$$

$$
U= H^{-1} (0,T ; L^2(-B,B)) + L^2 (0,T ; H^{-4}(-B,B)) +
L^{\frac{4}{3+\delta}}(0,T;L^{\frac{4}{3+\delta}}(-B,B)),
$$
and
$$
p=\frac{4}{3+\delta}.
$$

The proof is completed.

\subsection{Motivation and explanation of the main difficulty}
%\begin{remark}
Note that inclusions \eqref{limitations} can be obtained also for
$\delta=1$ with $\|u_0\|<1/2.$ Using embedding machinery and
interpolation theory for anisotropic spaces, one could pass to the
limit as $\eps \to 0$ in nonlinear term, as well. Indeed, let
$\delta=1.$ Multiplying $A^{\eps}u_{\eps}=0$ by $2(1+x)u_{\eps}$ and
integrating over $\Omega,$ we have

\begin{eqnarray*}\label{estimate 1.1}
 \frac{d}{dt}\left((1+x), u^2\right) (t) + \|\nabla
 u\|^2(t)+2\|u_x\|^2(t)+
 (1-2\eps)\int_{-B}^B u^2_x (0,y,t) \, dy\\
 \le \|u\|^2(t)+2\|u\|^4_{L^4(\Omega)}
 \le \|u\|^2(t)+2\|\nabla u\|^2(t)\|u\|^2(t).
\end{eqnarray*}
Bearing in mind that $\|u\|(t)\le \|u_0\|(t)<1/2$ and integrating in
$t>0,$ Gronwall's lemma gives $$u\in
L^{\infty}\left(0,T;L^2(\Omega)\right) \cap
L^2\left(0,T;H^1_0(\Omega)\right)$$ with both estimates independent
of $\eps<1/4.$

%________________________________________________________________

Now we observe that
\begin{eqnarray*}
\int_0^T \int_{\Omega} |u^3|^{\frac{4}{3}}  dx dt \leq C \|u_0\|^2
\| \nabla u\|^2_{L^2_T L^2_{xy}}
\end{eqnarray*}
%====================================
and by estimate above this implies $ u^3 \in L^{\frac{4}{3}}(Q_T). $
%====================================
Since $L^{\frac{4}{3}}(\Omega) \hookrightarrow H^{-1}(\Omega),$ we
conclude that
\begin{equation*}
u^2u_x = \frac{1}{3}\partial_x(u^3) \in L^{\frac{4}{3}}(0,T;
H^{-2}(\Omega))
\end{equation*}
%====================================
%Due to $ u_t= - \Delta u_x - u_x - u^2u_x $
%====================================
whence
\begin{equation*}
u_t \in L^{\frac{4}{3}}(0,T; H^{-2}(\Omega))
\end{equation*}
%_________________________________________________________________
and passage to the limit as $\eps \to 0$ in nonlinear term can be
justified as above.

It is difficult, however, to obtain explicit estimates like
\eqref{estimativa 3.3} with $m>1$ for $\delta=1.$
%which turns
%unclear how to obtain the estimates of higher order (and, partially,
%how to prove the uniqueness).
In fact, let $r,s \ge 1.$ We are going
to determine conditions upon $r$ and $s$ such that $u^2 u_x $ lies
in $L^r\left((0,T;L^s(\Omega)\right).$ Consider $p,q>1$ with $1/p
+1/q = 1.$ Then
\begin{eqnarray}\label{3.12}
\|u^2u_x \|_{L^r_T L^s_{xy}}^r &=& \int_0^T \left( \int_{\Omega} u^{2s} u_x^s \, d\Omega \right) ^{\frac{r}{s}} dt
\nonumber \\
%====
&\leq & \int_0^T
\|u\|^{2r}_{L^{2sp}_{xy}}(t)\|u_x\|^{r}_{L^{sq}_{xy}}(t) \, dt .
\end{eqnarray}
%=====
By Nirenberg's inequality with $
%\begin{equation}
\alpha = \frac{sp-1}{sp} $
%\end{equation}
one has
\begin{equation*}
\|u\|^{2r}_{L^{2sp}_{xy}}(t) \leq C \|\nabla u\|^{2r \alpha}
\|u\|^{2r(1-\alpha)}.
\end{equation*}
Supposing $sq\leq 2,$ estimate \eqref{3.12} reads
\begin{eqnarray*}
\|u^2u_x \|_{L^r_T L^s_{xy}}^r &\leq & C \|u \|_{L^{\infty}_T L^2_{xy}}^{2r(1-\alpha)} \int \|\nabla u\|^{2r \alpha} \|u_x\|^{r}(t) \, dt \nonumber \\
&\leq & C \|u \|_{L^{\infty}_T L^2_{xy}}^{2r(1-\alpha)} C \|\nabla u
\|_{L^{r(2\alpha +1)}_T L^2_{xy}}^{r(2\alpha +1)}.
\end{eqnarray*}
In order to gain $r(2\alpha +1) =2,$ it should be $\alpha = 1/r
-1/2.$ Therefore,
%\begin{eqnarray} \frac{1}{r} - \frac{1}{2} = 1 -
%\frac{1}{sp}  \,\,\, \Rightarrow \,\,\,
$
\frac{1}{sp}=  \frac{3}{2} - \frac{1}{r}, $
%\end{eqnarray}
which implies
\begin{eqnarray*}
sq = \frac{2rs}{2(r+s) - 3rs}.
\end{eqnarray*}
Since $sq \leq 2,$ it follows that
%\begin{eqnarray}
$\frac{2rs}{2(r+s) - 3rs} \leq 2$ which means  $sr \leq \frac{
r+s}{2}.$
%\end{eqnarray}
Observe that for $r,s>1$ this condition does not hold. The only
possibility thus reads $r=s=1,$ i.e., $u^2 u_x  \in
L^1\left((0,T;L^1(\Omega)\right).$

The space $(L^1_t;L^1_{xy})$ is known to be difficult to deal with.
For example, it is not clear even whether the condition
$u_x(L,y,t)=0$ being satisfied. We leave it here only to illustrate
a challenge appearing in the critical case.
%\end{remark}

\section{Local result for critical case}\label{local existence}

Consider the following Cauchy problem in abstract form:
\begin{equation}\label{Eq sem. grp.}
 \left\{
\begin{array}{c}
%-------------------------------
u_t + Au= f, \;\; \\ %\text{em} \;\;  \mathrm{Q}_T = \Omega\times (0,T), \;\; ,\\
%-------------------------------
%\hspace{-0.5cm}u=0\;\;\text{em}\;\;\partial \Omega , \;\; u_x(L,y,t)=0 , \\
%-------------------------------
\hspace{0cm}u(0)=u_0 ,\;\;%\text{em}\;\;\Omega ;
%-------------------------------
\end{array}
\right.
\end{equation}
where $f\in L^1(0,T;L^2(\Omega))$ and $A: L^2(\Omega)\to
L^2(\Omega)$ defined as $A\equiv
\partial_x + \Delta\partial_x $ with the domain
$$
D(A) = \{ u \in L^2(\Omega ) \, ; \, \Delta u_x + u_x \in
L^2(\Omega) \, \text{with}\, u|_{\partial \Omega}=0 \ \text{ and }\
u_x(L,y,t) = 0, \ t\in (0,T) \},
$$
endowed with its natural Hilbert norm $ \|u\|_{D(A)}(t) = \left(
\|u\|^2_{L^2(\Omega)}(t) + \| \Delta u_x + u_x \|^2_{L^2(\Omega)}(t)
\right)^{1/2} $ for all $t\in (0,T)$.
%\begin{proposition}\label{proposition1} $D(A) \hookrightarrow  H^1_0(\Omega)\cap H^2(\Omega)$ com imersao continua.
%\end{proposition}
%Prova em \cite{temam2}.
\begin{proposition}\label{res. 1 teman} Assume $u_0\in D(A)$ and $f\in L^1_{loc}(\R^+;L^2(\Omega))$
with $f_t\in L^1_{loc}(\R^+;L^2(\Omega))$. Then problem \eqref{Eq
sem. grp.} possesses the unique solution $u(t)$ such that
\begin{equation}\label{res 2}
u\in C([0,T]; D(A)) , \; u_t \in L^{\infty}(0,T; L^2(\Omega)) \;
T>0.
\end{equation}
Moreover, if $u_0\in L^2(\Omega)$ and $f \in
L^1_{loc}(\R^+;L^2(\Omega)),$ then \eqref{Eq sem. grp.} possesses a
unique (mild) solution $u \in C([0,T]; L^2(\Omega)) $ given by
\begin{equation}\label{res 3}
u(t)=S(t)u_0 + \int_0^t S(t-s)f(s)\, ds.
\end{equation}
\end{proposition}
\begin{corollary}\label{res. 2 teman} Under the hypothesys of Proposition \ref{res. 1 teman}, the solution
$u$ in (\ref{res 2}) satisfies
\begin{equation}\label{res 2.1}
u\in L^{\infty}((0,T);H^1_0(\Omega)\cap H^2(\Omega)),
\end{equation}
\end{corollary}
For the proof, see \cite{temam2}.

Furthermore, one can get (see \cite{kato}, for instance) the
estimate for strong solution (\ref{res 2}):
\begin{equation}\label{estimativa temam/kato 1}
\|u_t \| (t) \leq \|Au_0\| + \|f\|(0) + \|f_t\|_{L^1_tL^2_{xy}},
\end{equation}
and
\begin{equation}\label{estimativa temam/kato 2}
\left\| A u \right\| (t) \leq \|  u_t \|(t)+ \|f\|(t).
\end{equation}
Since $D(A){ \hookrightarrow }  H^1_0(\Omega)\cap H^2(\Omega)$
compactly (see \cite{temam2} for instance), we have the estimate
\begin{eqnarray}\label{estimativa temam/kato 4}
\|u\|_{L^{\infty}0,T; H^1_0\cap H^2 (\Omega)}(t) \leq C
\big(\|u\|_{L^{\infty}_t L^2_{xy}}
+  \|Au_0\| + \|f\|(0) + \|f_t\|_{L^1_tL^2_{xy}}  +  \|f\| _{L^{\infty}_t L^2_{xy}} \big). \\
\end{eqnarray}
%=======================
where $C$ depends only on $\Omega $. Next, we define
$$
Y_T = \{f \in L^{1}\big( 0,T; L^2(\Omega) \big) \; \text{such that}
\; f_t \in L^{1} \big( 0,T; L^2(\Omega)   \big)     \}
$$
with the norm
$$
\|f\|_{Y_T}= \|f\|_{ L^{1}_t L^2_{xy} }  + \|f_t\|_{L^{1}_t L^2_{xy}
}.
$$
\begin{remark}\label{prop traÃ§o f} If  $f\in Y_T,$ then $ f \in C([0,T];L^2(\Omega))$,
with the constant $C_T$ from $\|f\|_{C_t L^2_{xy}} \leq C_T
\|f\|_{Y_T} $ which is proportional to $T$ and its positive powers
\cite{evans}.
\end{remark}
Consider $X_T^0 = L^{\infty}\big(0,T; H^1_0(\Omega)  \cap
H^2(\Omega)  \big) $ and define the Banach space
\begin{eqnarray}\label{esp XT}
X_T=\{u \in X_T^0:\ u_t \in L^{\infty}\big(0,T; L^2(\Omega)\big)
\;\;\text{and}\;\; \nabla u_t \in L^{2}\big(0,T; L^2(\Omega)\big)
\} .
\end{eqnarray}
%-------------------------------
with the norm
\begin{eqnarray}\label{morma esp XT}
\| u \|_{X_T} = \| u \|_{L^{\infty}_T H^1_0 \cap H^2_{xy} }
+ \| u_t \|_{L^{\infty}_T L^2_{xy}}
+\| \nabla u_t \|_{L^{2}_T L^2_{xy}} . \\
\end{eqnarray}
%-------------------------------

\begin{thm}\label{Teo solu local mZK} Let $u_0 \in D(A)$. Then there exists $T>0$
such that IBVP (\ref{2.1})-(\ref{2.4}) possesses a unique
%(weak/mild)
solution in $X_T$.
\end{thm}
The proof of the Theorem consists in three lemmas below.
\begin{lemma}\label{lemma1} The function
%\begin{eqnarray*}
$Y_T
%&
\longrightarrow
%&
X_T;\ \
%\\
%-------------------------------
f
%&
\mapsto
%&
\int_0^t S(t-s)f(s) ds$
%\end{eqnarray*}
is well defined and
continuous.
\end{lemma}
%-------------------------------
%-------------------------------
%\noindent{{\sl Proof:}}
For the proof, note that this function maps $f$ to the solution of
homogeneous linear problem with zero initial datum. Estimates
(\ref{estimativa temam/kato 1}) and (\ref{estimativa temam/kato 4})
then give
\begin{equation}\label{estimativa ut 0}
\| u \|_{L^{\infty}_T H^1_0 \cap H^2_{xy} } + \| u_t \|_{L^{\infty}_T L^2_{xy}} \leq C \|f\|_{Y_T},
\end{equation}
where $C$ is as above. Thus, it rests to estimate the term $\|
\nabla u_t \|_{L^{2}_T L^2_{xy}}$ in (\ref{morma esp XT}).

Differentiate the equation in (\ref{Eq sem. grp.}) with respect to
$t,$ multiply it by $(1+x) u_{t}$ and integrate the outcome over
$\Omega.$ The result reads
\begin{equation}\label{estimativa ut}
\dfrac{d}{dt} \left( (1+x),u_t^2\right) (t) + \|\nabla u_t \|^2(t) +
2\| u_{xt} \|^2+\int_{-B}^B u_{xt}^2(0,y,t)\,dy = \|u_t\|^2(t) +
2\int_{\Omega}(1+x) f_t  u_{t} \, d\Omega.
\end{equation}
H\"older's inequality and (\ref{estimativa temam/kato 1}) imply
\begin{eqnarray}\label{estimativa ut 1.2}
\int_0^T  \|\nabla u_t \|^2(t)\, dt &\leq &  T\big(\|f\|(0) +
\|f_t\|_{L^1_T L^2_{xy}}  \big)^2 \nonumber \\
&+& 2(1+L)\big(\|f\|(0) + \|f_t\|_{L^1_T L^2_{xy}} \big)
\|f_t\|_{L^1_T L^2_{xy}} +\left( (1+x),u_t^2\right) (0) .
%==============================
\end{eqnarray}
%-------------------------------
%==============================
Using the equation from (\ref{Eq sem. grp.}) and taking in mind that
$u_0 \equiv 0 $, we get
\begin{eqnarray}\label{estimativa ut 1.3}
u_t (x,y,0) = f(x,y,0) - Au_0= f(x,y,0)
\end{eqnarray}
%=======================
Inserting (\ref{estimativa ut 1.3}) into (\ref{estimativa ut 1.2})
provides
\begin{eqnarray}\label{estimativa ut 1.7}
\|\nabla u_t \|^2_{L^2_TL^2_{xy}}  \leq  \Big(4T K_T^2 + 4K_T(1+L) + K_T^2(1+L)  \Big)\|f\|^2_{Y_T},
\end{eqnarray}
%==============================
where $K_T=\max \{1,C_T\}$. Therefore, estimates (\ref{estimativa ut
0}) and (\ref{estimativa ut 1.7}) read
\begin{eqnarray}\label{estimativa ut 1.9}
\|u\|_{X_T} \leq K \|f\|_{Y_T}.
\end{eqnarray}
%==============================
\begin{lemma} The function
$$
D(A)   \longrightarrow X_T ;\ u_0  \mapsto  S(t)u_0
$$
is well defined and continuous.
\end{lemma}
%-------------------------------
The proof follows the same steps as Lemma \ref{lemma1}, taking into
account that now $f\equiv 0$. The resulting estimate is
\begin{eqnarray}\label{estimativa linear XT}
 \|u\|_{X_T} & \leq &  M \|u_0 \|_{D(A)},
\end{eqnarray}
%-------------------------------
where $M$ is given by
\begin{eqnarray}\label{estimativa linear M}
M = 2C + 1 + \sqrt{1+L+T},
\end{eqnarray}
and $C$ (which depends only on $\Omega$) is defined by continuous
immersion $D(A) \hookrightarrow  H^1_0(\Omega)\cap H^2(\Omega).$
%-----------------------------------------------------------------

\begin{lemma}\label{contraÃ§Ã£o}
Given $R>0$, consider the closed ball $ B_R = \{u \in X_T  ;
\|u\|_{X_T} \leq R \}.$ Then the operator
\begin{eqnarray}\label{contraÃ§ao 0}
\Phi : B_R  \longrightarrow  X_T ;\
%-------------------------------
v \mapsto S(t)u_0 -\int_0^t S(t-s)v^2 v_x (s) \,ds \nonumber
\end{eqnarray}
is the contraction.
\end{lemma}
%-------------------------------
Fix $R>0$ and $u,v \in B_R.$ We have
\begin{eqnarray*}
\Phi (v) - \Phi (u)= \int_0^t S(t-s)[u^2u_x - v^2 v_x] (s) \,ds
\nonumber
\end{eqnarray*}
%-------------------------------
so that (\ref{estimativa ut 1.9}) implies
\begin{eqnarray}\label{contraÃ§ao 0.1}
\| \Phi (u) - \Phi (v) \|_{X_T} \leq K \|u^2u_x - v^2v_x\|_{Y_T}.
\end{eqnarray}
%-------------------------------
We study the right-hand norm in detail:
\begin{eqnarray}\label{contaÃ§ao 1}
\|u^2u_x - v^2v_x\|_{Y_T} &=& \|u^2u_x - v^2v_x\|_{L^{1}_T L^2_{xy}}
+\left\|\big(u^2u_x\big)_t - \big(v^2v_x\big)_t\right\|_{L^{1}_Y L^2_{xy}} \nonumber \\
&=& I + J.
\end{eqnarray}
%=============================
First, we write
\begin{eqnarray}\label{contaÃ§ao 2}
I &=& \left\|(u^2 - v^2)u_x\right\|_{L^{1}_T L^2_{xy}} +\left\| v^2(u_x-v_x)\right\|_{L^{1}_T L^2_{xy}}   \nonumber \\
&=& I_1 + I_2 .
\end{eqnarray}
For the integral $I_1$ one has
\begin{eqnarray}\label{contaÃ§ao 2.1}
I_1 \leq  \int_0^T  \| u-v\|_{L^{6}(\Omega)} \|
u+v\|_{L^{6}(\Omega)} \|u_x \|_{L^{6}(\Omega)} dt.
\end{eqnarray}
%============================
Nirenberg's inequality gives
\begin{eqnarray}\label{contaÃ§ao 2.4}
I_1  &\leq & T C_{\Omega} \| \nabla (u  + v)\|^{\frac{2}{3}}_{L^{\infty}_T L^2_{xy}} \|u
+ v\|^{\frac{1}{3}}_{L^{\infty}_T L^2_{xy}} \| \nabla u_x\|^{\frac{2}{3}}_{L^{\infty}_T L^2_{xy}}
\|u_x\|^{\frac{1}{3}}_{L^{\infty}_T L^2_{xy}} \| \nabla (u - v)\|^{\frac{2}{3}}_{L^{\infty}_T L^2_{xy}}
\|u - v\|^{\frac{1}{3}}_{L^{\infty}_T L^2_{xy}} \nonumber \\
%=============================
&= & T C_{\Omega}  D^{\frac{2}{3}} \| u  + v\|_{X_T} \|u\|_{X_T} \|
u - v \|_{X_T} ,
\end{eqnarray}
%=============================
where $D$ is the Poincare's constant from $\|w\|\leq D \| \nabla
w\|.$ Since $u$ and $v$ lie in $B_R,$ we conclude
\begin{eqnarray}\label{contaÃ§ao 2.5}
I_1 \leq  TK_0R^2 \| u - v \|_{X_T}.
\end{eqnarray}
%=============================
The integral $I_2$ can be treated in the similar way as $I_1$. It
rests to estimate the integral $J$.
\begin{eqnarray}\label{contaÃ§ao 3}
J \leq  \| 2u u_t(u_x - v_x)  \|_{L^1_T L^2_{xy}} + \| u^2 (u_{xt} - v_{xt}) \|_{L^1_T L^2_{xy}} + \|2v_x u(u_t - v_t)  \|_{L^1_T L^2_{xy}} \nonumber \\
%=============================
+ \|2v_x v_t (u-v)  \|_{L^1_T L^2_{xy}} +\| v_{xt}(u-v)(u+v) \|_{L^1_T L^2_{xy}} \nonumber \\
%=============================
= J_1 + J_2 + J_3 +J_4 + J_5.
\end{eqnarray}
%=============================
For $J_1$ we have
\begin{eqnarray}\label{contaÃ§ao 3.1}
J_1 & \leq \int_0^T  \| u\|_{L^{6}(\Omega)} \| u_t\|_{L^{6}(\Omega)} \|u_x - v_x \|_{L^{6}(\Omega)} \, dt .
\end{eqnarray}
%=============================
Niremberg's inequality implies
\begin{eqnarray}\label{contaÃ§ao 3.2}
J_1  \leq   C_{\Omega} \| \nabla u\|^{\frac{2}{3}}_{L^{\infty}_T L^2_{xy}} \|u\|^{\frac{1}{3}}_{L^{\infty}_T L^2_{xy}} \| \nabla (u_x - v_x)\|^{\frac{2}{3}}_{L^{\infty}_T L^2_{xy}} \|u_x-v_x\|^{\frac{1}{3}}_{L^{\infty}_T L^2_{xy}}\|u_t\|^{\frac{1}{3}}_{L^{\infty}_T L^2_{xy}} \nonumber \\
\leq T^{\frac{2}{3}} K_2 R^2  \| u - v \|_{X_T}.
\end{eqnarray}
%=============================
The integrals $J_3$ and $J_4$ are analogous to $J_1$. To get bound
for $J_5$ we observe that
\begin{eqnarray}\label{contaÃ§ao 3.6}
J_5  &=& \int_0^T \left(\int_{\Omega} v_{xt}^2 (u -v)^2 (u+v)^2 \, d\Omega\right)^{\frac{1}{2}} dt \nonumber \\
%==============================
& \leq & \int_0^T \left( \sup (u-v)^2 \right)^{\frac{1}{2}} \left( \sup (u+ v)^2 \right)^{\frac{1}{2}} \|v_{xt}\|(t) \, dt \nonumber \\
%==============================
&\leq & \int_0^T \big( \|u-v\|^2_{H^1_{xy}}(t) + \|u_{xy} - v_{xy}\|^2 (t) \big)^{\frac{1}{2}} \big( \|u+v\|^2_{H^1_{xy}}(t) + \|u_{xy} + v_{xy}\|^2 (t) \big)^{\frac{1}{2}}   \| v_{xt}\|(t) \, dt \nonumber \\
%=============================
&\leq & \big( \|u-v\|_{L^{\infty}_T H^1_{xy}} + \|u_{xy}- v_{xy}\|_{L^{\infty}_T L^2_{xy}} \big)\big( \|u+v\|_{L^{\infty}_T H^1_{xy}} + \|u_{xy}+ v_{xy}\|_{L^{\infty}_T L^2_{xy}} \big)\| v_{xt}\|_{L^{1}_T L^2_{xy}} \nonumber \\
%=============================
& \leq & 4 T^{\frac{1}{2}} \| v\|_{X_T}\|u+v\|_{X_T} \|u-v\|_{X_T} \nonumber \\
%=============================
& \leq & 8 T^{\frac{1}{2}} R^2 \|u-v\|_{X_T}.
\end{eqnarray}
%=============================
The integral $J_2$ follows like $J_5$. Thus,
\begin{eqnarray}\label{contaÃ§ao 3.6}
\|u^2u_x - v^2v_x\|_{Y_T} \leq K K^* T^{\frac{1}{2}}R^2 \|u-v\|_{X_T}.
\end{eqnarray}
%=============================
Finally, choosing $T>0$ such that $K K^* T^{\frac{1}{2}}R^2<1 ,$ we
conclude that $\Phi $ is a contraction map.

Lemma \ref{contraÃ§Ã£o} is proved.

Let $u \in B_R$. If $R= 2 M\|u_0\|_{D(A)},$ then estimates
(\ref{estimativa linear XT}) and (\ref{contaÃ§ao 3.6}) with $v
\equiv 0$ assure
\begin{eqnarray}\label{contaÃ§ao 3.7}
\|u\|_{X_T} & \leq & \|S(t)u_0\|_{X_T} +\| \int_0^t S(t-s)u^2u_x \, ds \|_{X_T} \nonumber \\
%=============================
& \leq & M\|u_0\|_{D(A)} + K K^* T^{\frac{1}{2}}R^2 \|u\|_{X_T} \nonumber \\
%=============================
&\leq & \frac{R}{2} + K K^* T^{\frac{1}{2}}R^3.
\end{eqnarray}
%=============================
Setting $T>0$ such that $K K^* T^{\frac{1}{2}}R^3 < \frac{R}{2}, $
one get
\begin{eqnarray}\label{contaÃ§ao 3.8}
\|u\|_{X_T} \leq R.
\end{eqnarray}
Choose $T>0$ such that $K K^* T^{\frac{1}{2}}R^2<1 $ and $ K K^*
T^{\frac{1}{2}}R^3 < \frac{R}{2} .$ Then $\Phi$ is the contraction
from the ball $B_R$ into itself. Therefore, the Banach fixed point
theorem assures the existence of a unique element $u\in B_R$ such
that $\Phi (u) = u.$

This completes the proof of Theorem \ref{Teo solu local mZK}.

\section{Decay}
%-------------------------------

\begin{thm}
Let $B,L>0$ satisfy $$\pi^2
\left[\frac{3}{L^2}+\frac{1}{4B^2}\right] - 1 :=2 A^2 >0\ \ \text{
and
 }\ \ \|u_0\|^2 < \dfrac{A^2}{2\pi^2\left(
\frac{1}{L^2}+\frac{1}{4B^2}\right)}.$$ If there exists solution
$$
u\in L^{\infty}\left(0,\infty;H^1_0(\Omega)\right) $$ to
\eqref{2.1}-\eqref{2.4}, then
\begin{equation}\label{decaimento 1}
%-------------------------------
\|u\|^2(t) \leq \left( 1+x, u^2 \right) (t) \leq
e^{-\left(\frac{A^2}{(1+L)}\right)t}\left( 1+x, u^2_0 \right).
%-------------------------------
\end{equation}
\end{thm}
%-------------------------------
%-------------------------------
%-------------------------------

%-------------------------------
%-------------------------------
To prove this result we will use
\begin{lemma}  \label{prop 1 decaimento} \text{\sl (V. A. Steklov)}
Let $L,B>0$ and $\omega \in H_0^1(\Omega)$. Then
\begin{equation}\label{decaimento 1.4}
\int_0^L \int_{-B}^B \omega^2(x,y)dxdy \leq
\frac{4B^2}{\pi^2}\int_0^L \int_{-B}^B \omega_y^2(x,y)dxdy,
\end{equation}
%-------------------------------
%-------------------------------
and

\begin{equation}\label{decaimento 1.5}
\int_0^L \int_{-B}^B \omega^2(x,y)dxdy \leq
\frac{L^2}{\pi^2}\int_0^L \int_{-B}^B \omega_x^2(x,y)dxdy.
\end{equation}
\end{lemma}
See \cite{doronin} for the proof.
%-------------------------------
%-------------------------------
We start the proof of \eqref{decaimento 1}, multiplying \eqref{2.1}
by $u$ and integrating over $Q_t,$ which easily gives

\begin{equation}\label{decaimento 1.3}
\| u\| ^2 (t) \leq \|u_0\|^2.
\end{equation}

%-------------------------------
%-------------------------------
%-------------------------------

Multiplying \eqref{2.1} by $(1+x)u$ and integrating over $\Omega,$
we have

\begin{eqnarray}\label{decaimento 1.6}
\dfrac{d}{dt} \left( 1+x, u^2 \right)(t) + \int_{-B}^B u^2_x (0,y,t)
\, dy
+  \|\nabla u \|^2(t) + 2\|u_x\|^2(t) - \|u\|^2(t) \nonumber \\
%-------------------------------
=   -2\int_{\Omega}(1+ x) u (u^2u_x) \, d\Omega %\nonumber \\
%-------------------------------
=  \frac{1}{2}\int_{\Omega}u^4 \, d\Omega.%\nonumber \\
\end{eqnarray}
%-------------------------------
%-------------------------------
For the integral $I_1 = \frac{1}{2}\int_{\Omega}u^4 =
\frac{1}{2}\|u\|^4_{L^4(\Omega)}(t),$
%\begin{eqnarray}\label{decaimento 1.7}
%\end{eqnarray}
%-------------------------------
%-------------------------------
Nirenberg's inequality implies

\begin{eqnarray}\label{decaimento 1.8}
I_1 &\leq & \frac{1}{2} \big(  2^{\frac{1}{2}} \|\nabla u \|^{\frac{1}{2}}(t)\| u \|^{\frac{1}{2}} (t) \big)^4 \nonumber \\
%-------------------------------
&=& 2 \|\nabla u \|^{2}(t)\| u \|^{2}(t)%\nonumber \\
%-------------------------------
\leq  2 \|\nabla u \|^{2}(t)\| u_0\|^{2}(t).
\end{eqnarray}
%-------------------------------
%-------------------------------
Take
$$ I_2 = 3\|u_x\|^2(t)+ \|u_y\|^2(t).
$$
%-------------------------------
%-------------------------------
For all $\varepsilon>0$ we have
$$ I_2 = (3-\varepsilon)\|u_x\|^2(t)+ (1-\varepsilon)\|u_y\|^2(t) + \varepsilon\big( \|u_x\|^2(t)+ \|u_y\|^2(t) \big).
$$
%-------------------------------
%-------------------------------
Lemma \ref{prop 1 decaimento} jointly with (\ref{decaimento 1.6})
and (\ref{decaimento 1.8}) provides

\begin{eqnarray}\label{decaimento 1.9}
\dfrac{d}{dt} \left( 1+x, u^2 \right)(t)
+ \left[ \pi^2\left( \frac{3}{L^2}+\frac{1}{4B^2} \right) -1 - \varepsilon \pi^2\left( \frac{1}{L^2}+\frac{1}{4B^2} \right)  \right]\|u\|^2(t)  \nonumber \\
%-------------------------------
+ \left( \varepsilon - 2\|u_0\|^2 \right)\|\nabla u\|^2(t) \leq 0.
\end{eqnarray}
%-------------------------------
%-------------------------------
Define
$$
2 A^2:=\pi^2 \left[\frac{3}{L^2}+\frac{1}{4B^2}\right] - 1 >0,
%$$
\ \text{ and take }\
%$$
\varepsilon= \dfrac{A^2}{\pi^2\left(
\frac{1}{L^2}+\frac{1}{4B^2}\right)}.
$$
The result for \eqref{decaimento 1.9} reads
\begin{eqnarray}\label{decaimento 1.10}
\dfrac{d}{dt} \left( 1+x, u^2 \right)(t) + A^2\|u\|^2(t)+ \left(
\varepsilon - 2\|u_0\|^2 \right)\|\nabla u\|^2(t) \leq 0.
\end{eqnarray}
%-------------------------------
%-------------------------------
If $0\leq \varepsilon - 2\|u_0\|^2,$ then
\begin{eqnarray}\label{decaimento 1.11}
\dfrac{d}{dt} \left( 1+x, u^2 \right)(t) + \frac{A^2}{(1+L)}\left(
1+x, u^2 \right)(t) \leq 0,
\end{eqnarray}
%-------------------------------
%-------------------------------
and consequently
\begin{equation}\label{decaimento 1.12}
%-------------------------------
\|u\|^2(t) \leq \left( 1+x, u^2 \right) (t) \leq
e^{-\left(\frac{A^2}{(1+L)}\right)t}\left( 1+x, u^2_0 \right).
%-------------------------------
\end{equation}

The proof is completed.

%\section{Acknowledgments}
%We appreciate very much fruitful and constructive comments and
%suggestions of both of the Reviewers.

\end{document}